\documentclass{amsart}
\usepackage[utf8]{inputenc}
\usepackage{fancyhdr}
\usepackage{chngcntr}
\usepackage{multicol}
\makeindex

\usepackage{mathtools}

\usepackage{enumitem}
\usepackage[colorlinks, linkcolor=blue]{hyperref}

\usepackage[ngerman,american]{babel}
\usepackage[T1]{fontenc}
\usepackage{upgreek}
\usepackage{tikz}
\usepackage{float}
\usetikzlibrary{matrix}
\usepackage[arrow, matrix, curve]{xy}

\usepackage{epstopdf}

\usepackage{dsfont}
\newcommand{\bo}[3]{\begin{#2}[#1]#3\end{#2}}
\newcommand{\bs}[2]{\begin{#1}#2\end{#1}}
\newcommand{\bd}[3]{\begin{#1}{#2}#3\end{#1}}
\newcommand{\be}[3][]{\ifthenelse{\equal{#1}{}}{\begin{#2}#3\end{#2}}{\begin{#2}[#1]#3\end{#2}}} 
\newcommand{\test}[3][]{\ifthenelse{\equal{#1}{}}{#1}{#1}} 

\newcommand{\ft}[1]{\text{#1}} 
\newcommand{\ftn}[1]{\textnormal{#1}} 

\newif\ifSIMPLE
\newcommand{\ftBranch}[2]{\ifSIMPLE {#1}\else{#2}\fi}

\usepackage{xspace}

\newcommand{\rhs}{right-hand side\xspace}
\newcommand{\wrt}{with respect to\xspace}

\newcommand{\dfn}{\coloneqq}

\newcommand{\mb}[1]{\left(#1\right)} 
\newcommand{\mba}[1]{\left\langle#1\right\rangle} 
\newcommand{\mbs}[1]{\left\{#1\right\}} 
\newcommand{\mbro}[1]{\left[#1\right)} 
\newcommand{\mbm}[1]{\left|#1\right|} 

\newcommand{\mnp}[2][2]{\left\|#2\right\|_{#1}} 

\usepackage{amsthm}

\newcommand{\se}[2][]{\ifthenelse{\equal{#1}{}}{Section \ref{s:#2}}{Section \ref{s:#2:#1}}}

\newcommand{\fg}[1]{Figure \ref{f:#1}}

\newcommand{\rem}[2][]{\ifthenelse{\equal{#1}{}}{Remark \ref{r:#2}}{\tm{#2} \ref{r:#2:#1}}}
\newcommand{\tm}[2][]{\ifthenelse{\equal{#1}{}}{Theorem \ref{t:#2}}{\tm{#2} \ref{t:#2:#1}}}
\newcommand{\lem}[2][]{\ifthenelse{\equal{#1}{}}{Lemma \ref{l:#2}}{\lem{#2} \ref{l:#2:#1}}}
\newcommand{\co}[2][]{\ifthenelse{\equal{#1}{}}{Corollary \ref{c:#2}}{\co{#2} \ref{c:#2:#1}}}

\newcommand{\eq}[2][]{\ifthenelse{\equal{#1}{}}{Equation \ref{eq:#2}}{\eq{#2} \ref{eq:#2:#1}}}
\newcommand{\ex}[2][]{\ifthenelse{\equal{#1}{}}{Example \ref{ex:#2}}{\ex{#2} \ref{ex:#2:#1}}}
\newcommand{\de}[2][]{\ifthenelse{\equal{#1}{}}{Definition \ref{d:#2}}{\de{#2} \ref{d:#2:#1}}}
\newcommand{\pp}[2][]{\ifthenelse{\equal{#1}{}}{Proposition \ref{p:#2}}{\pp{#2} \ref{p:#2:#1}}}


\DeclareMathOperator{\iunit}{\ftn{i}}

\DeclareMathOperator{\id}{\mathds{1}}

\DeclareMathOperator{\tr}{Tr}
\DeclareMathOperator{\spec}{sp}

\newcommand{\SPA}[1][and]{\quad\ft{#1}\quad}

\newcommand{\DIFF}[2][x]{\ftn{D}_{#1}#2}

\newcommand{\FUNC}[3][T]{I_{#2,#3}\mb{#1}}

\newcommand{\SL}[2][\mathbb{R}]{\ftn{SL}_{#2}\mb{#1}}

\newcommand{\MAT}[2][\mathbb{R}]{\ftn{Mat}_{#2}\mb{#1}}

\newcommand{\SR}{\ft{ }\middle|\ft{ }}

\newcommand{\MC}[1]{\mathcal{#1}}
\newcommand{\MB}[1]{\mathbb{#1}}

\newcommand{\VS}[2][d]{\ifthenelse{\equal{#1}{1}}{#2}{\ifSIMPLE{#2^2}\else{#2^{#1}}\fi}}

\newcommand{\PERT}[1][B]{\ifSIMPLE{\MC{#1}_{r}}\else{\MC{#1}_{d,r}}\fi}
\newcommand{\TORUS}[1][d]{\ifthenelse{\equal{#1}{1}}{\MB{T}}{\ifSIMPLE{\MB{T}^2}\else{\MB{T}^{#1}}\fi}}
\newcommand{\STRIP}[2][d]{\ifSIMPLE{\MB{A}_{#2}}\else{\MB{A}_{#2,#1}}\fi}
\newcommand{\ANORM}[1]{\ifSIMPLE{\mnp[{}]{#1}}\else{\mbm{#1}_{M}}\fi}

\newcommand{\MNORM}[1]{\ifSIMPLE{\mnp[{}]{#1}}\else{\mnp[1]{#1}}\fi}

\usepackage[backend=bibtex8, style=ieee, doi=false, isbn=false, urldate=comp, dateabbrev=false, sorting=nty]{biblatex}
\addbibresource{bibliography.bib}
\AtEveryBibitem
{
	\clearlist{address}
	\clearlist{language}
	\clearfield{month}
}

\setlist[enumerate,1]{label=(\roman*),font=\normalfont}

\setcounter{section}{0}

\setlength\parindent{0pt}	
\SIMPLEtrue


%
	
\bs{document}
{
	\title[Generic non-trivial resonances for Anosov diffeomorphisms]{Generic non-trivial resonances for Anosov diffeomorphisms}

	\author{Alexander Adam}
	\address{Sorbonne Universit\'es, UPMC Univ. Paris 6, CNRS, Institut de Math\'ematiques de Jussieu,
		4, Place Jussieu, 75005 Paris, France}
	\email{alexander.adam@imj-prg.fr}
	\date{\today}
	
	\thanks{I thank Frédéric Naud for sharing his ideas \cite{naud_2015} which gave rise to this work and Paolo Giulietti and Carlangelo Liverani \cite{giul_liv_2016} for pointing out the relevance of non-trivial resonances. Special thanks to Viviane Baladi for all her valuable suggestions and our discussions.}
	\bs{abstract}
	{
		We study real analytic perturbations of hyperbolic linear automorphisms on the $2$-torus. The Koopman and the transfer operator are nuclear of order $0$ when acting on a suitable Hilbert space. We show the generic existence of non-trivial Ruelle resonances for both operators. We prove that some of the perturbations preserve the volume and some of them do not. 
	}
	\maketitle

	\noindent
	\section{Introduction}
Let $T\colon\TORUS\rightarrow\TORUS$ be a real analytic Anosov diffeomorphism. We define the Ruelle resonances of $T$ to be the zeroes of the (holomorphically continued in $z\in\MB{C}$) dynamical determinant
\bs{align}
{
	\label{e:DYNDET}
	d_T(z)\dfn \exp -\sum_{n=1}^\infty \frac {z^n}n\sum_{T^n\mb{x}=x}\mbm{\det \mb{\id-\DIFF[x]{T^n}} }^{-1}\ft{.}
}
 It is well-known that $1$ is the only resonance if $T$ is a hyperbolic linear toral automorphism $M$. A subset of the Banach space of $\TORUS$-preserving maps, holomorphic and uniformly bounded on some annulus, is called generic if it is open and dense. We show in \tm{NTRS}, using an idea of Naud \cite{naud_2015}, that there is such a set $\MC{G}$ so that for all $\psi\in\MC{G}$, appropriately scaled, the Anosov diffeomorphism $M+\psi$ admits non-trivial Ruelle resonances. For this, we construct a Hilbert space of anisotropic generalized functions on which the transfer operator $\MC{L}_Tf\coloneqq(f/|\det \DIFF[{}] T|) \circ T^{-1}$
 is nuclear with its Fredholm determinant equal to $d_T$. Moreover, we prove that some of those
 generic perturbations preserve the volume while some do not.\\
The expanding case is easier and was initially
studied by Ruelle \cite{Ruelle1976}.
More recently, Bandtlow et. al \cite{bandtlow_spectral, slipantschuk_analytic_2013} calculated the resonances of real analytic expanding maps $T\colon S\rightarrow S$ on the unit circle $S$ explicitly for Blaschke
products. Their transfer operator acts on the Hardy space of holomorphic functions
on the annulus. (See also Keller and Rugh \cite{keller_eigenfunctions_2004} in the differentiable category.)\\
In the hyperbolic setting, Rugh proved the holomorphy of the dynamical determinant of real analytic Anosov diffeomorphisms on surfaces \cite{rugh_correlation_1992, Rugh_1996}. The idea was generalized by Fried to hyperbolic flows in all dimensions \cite{fried1995meromorphic}. The detailed study of anisotropic Banach spaces in the hyperbolic case
started with the pioneering work of \cite{blank_ruelleperronfrobenius_2002} (in the
differentiable setting) and is now a well established tool, see e.g. \cite{baladi_dynamical_2008} and \cite{ goue_liv_2006}.\\
Faure and Roy \cite{faure2006ruelle} later addressed real analytic perturbations of hyperbolic linear toral automorphisms on the two-dimensional torus, considering an anisotropic complex Hilbert space, which had already been briefly discussed by Fried \cite[Sect 8, I]{fried1995meromorphic}.\\
Our approach is based on this construction and strongly relies on an idea suggested by Naud \cite{naud_2015}. We  put  the transfer operator central
in our analysis. We introduce an anisotropic Hilbert space (\de{HS}) in \se{HILBERT}.\\
In \se{TRACEKOOP}, we \ftBranch{rephrase}{extend} a result from Faure and Roy \cite[Theorem 6]{faure2006ruelle} to show that the Koopman operator $\MC{K}_Tf\coloneqq f\circ T$ is nuclear of order $0$ when acting on our anisotropic Hilbert space.\\
In \se{SPECRUELLE}, we use this result and an idea of Naud \cite{naud_2015} to show that the Koopman operator admits non-trivial Ruelle resonances under a small generic perturbation of the dynamics.\\
In \se{SPECTRANSFER}, we consider the adjoint of the Koopman operator, which is just the transfer operator, acting on the dual Hilbert space and obtain our final results.\\\newline
Blaschke products were recently generalized to the hyperbolic setting
by Slipantschuk et al. \cite{Slipantschuk_2016} who calculate the entire spectrum of these
real analytic Anosov volume preserving diffeomorphisms explicitly. 

	\section{An anisotropic Hilbert space}
\label{s:HILBERT}
We denote\ftBranch{}{ for every $d\in\MB{N}$} the flat \ftBranch{$2$}{$d$}-torus by \ftBranch{$\MB{T}^2\dfn \MB{R}^2/\MB{Z}^2$}{$\MB{T}^d\dfn \MB{R}^d/\MB{Z}^d$}. We embed $\VS{\MB{T}}$ into the standard polyannulus in $\VS{\MB{C}}$ and set for each $r>0$ \[\STRIP{r}\dfn \VS{\MB{T}}+\iunit\VS{\mb{-r,r}}\ft{.}\]
We see $\STRIP{r}$ as a submanifold of $\VS{\MB{C}}$.
The Hilbert space $L_2\mb{\TORUS}$ is equipped with the canonical Lebesgue measure on $\TORUS$. This space admits an orthonormal Fourier basis given by
\bs{align}
{
	\label{e:FBL2}
	\varphi_n&\colon\TORUS\rightarrow\MB{C}\colon x\mapsto \exp\mb{\iunit 2\pi n^*x}\ft{,}\quad n\in\VS{\MB{Z}}\ft{,}
}
where $n^*$ is the canonical dual of $n$. We recall a construction from Faure and Roy \cite{faure2006ruelle} for a complex Hilbert space $\MC{H}_{A_{M,c}}$. This space also has been described briefly by Fried as an "ad hoc example" \cite[Sect. 8, I.]{fried1995meromorphic} of a generalized function space. The construction will be based on:
\bo{Hardy space $H_2\mb{\STRIP{r}}$}{definition}
{
	 For each $r>0$ and each holomorphic function $f\colon\STRIP{r}\rightarrow\MB{C}$, we define the norm 
	\[\mnp[H_2\mb{\STRIP{r}}]{f}\dfn \sup_{y\in\VS{\mb{-r,r}}}\mb{\int_{\TORUS}\mbm{f\mb{x+\iunit y}}^2\mathrm{d}x}^{\frac 12}\ft{.}\]
	Then we set
	\[H_2\mb{\STRIP{r}}\dfn\mbs{f\colon \STRIP{r}\rightarrow\MB{C} \SR f\ftn{ holomorphic, }\mnp[H_2\mb{\STRIP{r}}]{f}<\infty }\ft{.}\]
}
The space $H_2\mb{\STRIP{r}}$ is the \ftBranch{$2$}{$d$}-dimensional analogue of the Hardy space studied in \cite[p. 4]{sarason_1965}.
It admits a Fourier basis given by 
\[\vartheta^r_n\colon\STRIP{r}\rightarrow\MB{C}\colon x\mapsto \exp\mb{-2\pi r\MNORM{n}}\varphi_n\ft{,}\quad n\in\VS{\MB{Z}}\ft{,}\]
where $\MNORM{z}\coloneqq\mbm{z_1}+\mbm{z_2}$ for all $\mb{z_1,z_2}\eqqcolon z\in\VS{\MB{C}}$ and $z\in\TORUS$. With this choice of norm, the Fourier basis is orthonormal. Under the canonical isomorphism $L_2\mb{\TORUS}\cong L_2\mb{\TORUS}^*$, we have the isomorphism
\bs{equation}
{
	\label{e:CI}
	\mb{\vartheta_n^r}^*\cong \vartheta_n^{-r}\ft{.}
}
\ftBranch
{A matrix $M\in\SL[\MB{Z}]{2}$ is called hyperbolic if its eigenvalues lie not on the unit circle. We denote by $E^+_M$ the eigenspace of the eigenvalue of modulus $\lambda_M>1$ and by $E^-_M$ the eigenspace of the eigenvalue of modulus $\lambda_M^{-1}$.}
{A matrix $M\in\MAT[\MB{Z}]{d}$ is called hyperbolic if its eigenvalues are away from the unit circle. We let $E^u$ be the generalized eigenspace of $M$ corresponding to eigenvalues $\lambda$ with $\mbm{\lambda}>1$ and analogously $E^u$.} 
We decompose $y\in\VS{\MB{R}}$ uniquely as
\bs{align}{
	\label{e:DECOM}
	y=y_M^++y_M^-\SPA[with]y_M^+\in E^+_{M^*}\SPA[and]y_M^-\in E^-_{M^*}\ft{.}
	}
\ftBranch
{
	We have \bs{equation}
	{\label{e:CE}
		\ANORM{M^*y_M^+}=\lambda_M\ANORM{y_M^+}\SPA\ANORM{M^*y_M^-}=\lambda_M^{-1}\ANORM{y_M^-}\ft{.}}
}
{
	We let \[\mbm{\cdot}_M\colon\MB{R}^d\rightarrow\MB{R}_{\ge 0}\] be any norm for which there is a norm dependent constant $\lambda_M>1$ such that
	\bs{equation}
	{
		\label{e:CE}
		\mbm{M^*y_M^+}_M\ge\lambda_M\mbm{y_M^+}_M\SPA\mbm{M^*y_M^-}_M\le\lambda_M^{-1}\mbm{y_M^-}_M\ft{.}
	}
	For a construction of such a norm see e.g. the proof of \cite[Proposition 5.2.2]{brin_introduction_2002} and the constant $\lambda_M$ also \cite[p. 11]{barreira_nonuniform_2007}. We call $\mbm{\cdot}_M$ an adapted norm \wrt the adjoint $M^*$.
}

\bo{Scaling map $A_{M,c}$}{definition}
{
	\label{d:IM}
	Let $c>0$,\ftBranch{ and $M\in\SL[\MB{Z}]{2}$}{, $d\in\MB{N}$ and let $M\in\MAT[\MB{Z}]{d}$} be hyperbolic. For every $n\in\VS{\MB{Z}}$, we set, recalling (\ref{e:FBL2}), 
	\[A_{M,c}\varphi_n\dfn \exp\mb{-2\pi c\mb{\mnp[{}]{n_M^+}-\mnp[{}]{n_M^-}}}\varphi_n\ft{.}\]
}
\bo{Continuous embedding of $H_2\mb{\STRIP{r}}$}{lemma}
{
	\label{l:CE}
	Let $c>0$\ftBranch{ and let $M\in\SL[\MB{Z}]{2}$}{, $d\in\MB{N}$ and let $M\in\MAT[\MB{Z}]{d}$} be hyperbolic. Then the map $A_{M,c}$ can be extended by continuity to an injective linear map
		\[A_{M,c}\colon H_2\mb{\STRIP{c}}\rightarrow L_2\mb{\TORUS}\ft{,}\]
		bounded in operator norm by $1$.
}
\bs{proof}
{
	By \de{IM}, for each $f\in H_2\mb{\STRIP{c}}$ we have
	\bs{align*}
	{
		\mnp[L_2\mb{\TORUS}]{A_{M,c}f}^2&=\sum_{n\in\MB{Z}^2}\mbm{\varphi_n^*A_{M,c}f}^2=\sum_{n\in\MB{Z}^2}\exp\mb{-4\pi c\mb{\mnp[]{n^+_M}-\mnp[]{n^-_M}}}\mbm{\varphi_n^*f}^2\\
	&=\sum_{n\in\MB{Z}^2}\exp\mb{-4\pi c\mb{\mnp[]{n^+_M}-\mnp[]{n^-_M}+\mnp[1]{n}}}\mbm{{\vartheta_{n}^{c}}^*f}^2\ft{,}
	}
	where we used (\ref{e:CI}) in the last step. Using the triangle inequality, we find
	\[\mnp[]{n^+_M}-\mnp[]{n^-_M}+\MNORM{n}\ge 0\ft{.}\]
	Hence, it holds
	\[\sum_{n\in\MB{Z}^2}\exp\mb{-4\pi c\mb{\mnp[]{n^+_M}-\mnp[]{n^-_M}+\MNORM{n}}}\mbm{{\vartheta_{n}^{c}}^*f}^2\le \mnp[H_2\mb{\STRIP{c}}]{f}^2\ft{.}\]
	Injectivity follows since $A_{M,c}$ is invertible on the Fourier basis of $L_2\mb{\TORUS}$.
}
The image of $H_2\mb{\STRIP{c}}$ under $A_{M,c}$ is dense in $ L_2\mb{\TORUS}$ since it contains all Fourier polynomials.
\bo{Hilbert space $\MC{H}_{A_{M,c}}$}{definition}
{
	\label{d:HS}
	Let $c>0$\ftBranch{ and let $M\in\SL[\MB{Z}]{2}$}{, $d\in\MB{N}$ and let $M\in\MAT[\MB{Z}]{d}$} be hyperbolic. Let $A_{M,c}$ be the map given by \de{IM}. Then we set 
	\[\MC{H}_{A_{M,c}}\dfn \ft{closure of }H_2\mb{\STRIP{c}}\ft{ \wrt the norm }\mnp[L_2\mb{\TORUS}]{A_{M,c}\cdot}\ft{,}\]
	and extend $A_{M,c}$ by continuity to a linear map
	\[A_{M,c}\colon \MC{H}_{A_{M,c}}\rightarrow L_2\mb{\TORUS}\ft{.}\]
}
As a direct consequence of this construction, the scalar product on $\MC{H}_{A_{M,c}}$ satisfies
\[\mba{\cdot,\cdot}_{\MC{H}_{A_{M,c}}}\colon \MC{H}_{A_{M,c}}\times \MC{H}_{A_{M,c}}\rightarrow\MB{C}\colon \mb{f,g}\mapsto\mba{A_{M,c}f,A_{M,c}g}_{L_2\mb{\TORUS}}\ft{.}\]
An orthonormal Fourier basis of $\MC{H}_{A_{M,c}}$ is given by
\bs{equation}
{
	\label{e:FB}
	\varrho_n\dfn A_{M,c}^{-1}\varphi_n\ft{,}\quad n\in\VS{\MB{Z}}\ft{.}
}	
\bo{Dual space of $\MC{H}_{A_{M,c}}$}{lemma}
{
	\label{l:DS}
	Under the canonical isomorphism $L_2\mb{\TORUS}\cong L_2\mb{\TORUS}^*$, the dual space $\MC{H}_{A_{M,c}}^*$ is isomorphic to $A_{M,c}^2\MC{H}_{A_{M,c}}$.
}
\bs{proof}
{
	Under the canonical isomorphism $L_2\mb{\TORUS}\cong L_2\mb{\TORUS}^*$, we have for each $n_1$, $n_2\in\VS{\MB{Z}}$, using (\ref{e:FB}),
	\[\varphi_{n_1}^*\mb{\varphi_{n_2}}=\varphi_{n_1}^*\mb{A_{M,c}\varrho_{n_2}}=\mb{A_{M,c}\varphi_{n_1}}^*\mb{\varrho_{n_2}}=\mb{A_{M,c}^2\varrho_{n_1}}^*\mb{\varrho_{n_2}}\ft{.}\]
}

\bs{remark}
{
	\label{r:DS}
	By \lem{DS}, we associate to every linear functional $f^*\in\MC{H}_{A_{M,c}}^*$ a unique vector $f\in A_{M,c}^2\MC{H}_{A_{M,c}}$. Then, for every $g\in\MC{H}_{A_{M,c}}$, the product $fg$ is absolutely integrable \wrt the Lebesgue measure on $\TORUS$.
}
The decomposition in (\ref{e:DECOM}) defines two cones
\[C_{M}^+\dfn\mbs{y\in\VS{\MB{R}}\SR \ANORM{y_M^+}\ge \ANORM{y_M^-}}\SPA C_{M}^-\dfn\mbs{y\in\VS{\MB{R}}\SR \ANORM{y_M^+}\le \ANORM{y_M^-}}\ft{.}\]
\ftBranch{}{ If $E^-_M$ is trivial so is $C_{M}^-$. This is only the case if $M$ is expanding.}
\bs{example}
{
	\label{ex:MAPM}
	We let $M=\bs{pmatrix}{3&1\\2&1}$, then $\lambda_M=2+\sqrt{3}$. An eigenvector for $\lambda_M$ for $M^*$ is $\mb{1+\sqrt{3},1}$ and an eigenvector for $\lambda_M^{-1}$ is $\mb{1-\sqrt{3},1}$. The two subspaces $E_{M^*}^+$ and $E_{M^*}^-$ and the two cones $C_M^+$ and $C_M^-$ are shown in \fg{CONES}.
}
\bo{h!}{figure}
{
	\includegraphics{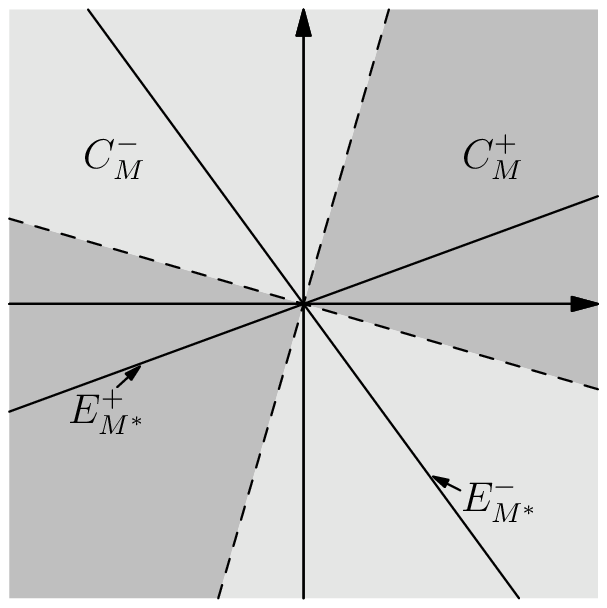}
	\caption{The map $M$ is from \ex{MAPM}. The dark gray area is the cone $C_M^+$ which contains the subspace $E_{M^*}^+$. The light gray area is the cone $C_M^-$ and contains $E_{M^*}^-$. A part $y\in\MB{R}^2$ belongs to the dashed lines if and only if $\ANORM{y_M^+}= \ANORM{y_M^-}$.}
	\label{f:CONES}
}
We set
\bs{align*}
{
	\MC{H}_{A_{M,c}}^+&\dfn\mbs{\sum_{n\in C_{M}^+\cap\VS{\MB{Z}}}\mba{\varrho_n,f}_{\MC{H}_{A_{M,c}}}\varrho_n\SR f\in\MC{H}_{A_{M,c}}}\ft{ and}\\
	\MC{H}_{A_{M,c}}^-&\dfn\mbs{\sum_{n\in C_{M}^-\cap\VS{\MB{Z}}}\mba{\varrho_n,f}_{\MC{H}_{A_{M,c}}}\varrho_n\SR f\in\MC{H}_{A_{M,c}}}\ft{.}
}
Hence, we have $\MC{H}_{A_{M,c}}=\MC{H}_{A_{M,c}}^++\MC{H}_{A_{M,c}}^-$. Comparing for each $n\in C_M^-$ the Fourier basis $\varrho_n$ with $\varphi_n$, it follows immediately that $\MC{H}_{A_{M,c}}^-\subset L_2\mb{\TORUS}$. For each $n\in C_M^+$, comparing the Fourier basis $\varrho_n$ with ${\vartheta_n^{c}}^*$, using (\ref{e:CI}), shows $\MC{H}_{A_{M,c}}^+\subset H_2\mb{\STRIP{c}}^*$. We conclude therefore that $\MC{H}_{A_{M,c}}$ contains linear functionals which do not belong to $L_2\mb{\TORUS}$. \ftBranch{}{In the case that $M$ is expanding we have $\MC{H}_{A_{M,c}}=\MC{H}_{A_{M,c}}^$ and by \lem{DS} it holds $\mb{\MC{H}_{A_{M,c}}}^*\cong A_{M,c}^2\MC{H}_{A_{M,c}}$. Now $A_{M,c}^2\MC{H}_{A_{M,c}}\subset L_2\mb{\TORUS}$ implying that there is no need for this being interpreted in terms of functionals.}
By construction, the space $\MC{H}_{A_{M,c}}$ is a rigged Hilbert space, i.e.: 
\bs{align}
{
	\label{e:RHS}
	H_2\mb{\STRIP{c}}\subset \MC{H}_{A_{M,c}}\subset H_2\mb{\STRIP{c}}^*\ft{.}
}
\bs{remark}
{
	We note that in the construction of $\MC{H}_{A_{M,c}}$, the expanding and contracting directions appear in the dual coordinates $n\in\VS{\MB{Z}}$ of the Fourier basis (\ref{e:FB}). This distinguishes $\MC{H}_{A_{M,c}}$ from the space of Rugh \cite{rugh_correlation_1992} where expanding and contracting coordinates are spatial. We observe
	\[n^*x=\mb{n_M^++n_M^-}^*\mb{x_{M^*}^++n_{M^*}^-}=\mb{n_M^+}^*x_{M^*}^++\mb{n_M^-}^*x_{M^*}^-\ft{.}\]
	Hence, we can rewrite (\ref{e:FB}) as
	\bs{align}
	{
		\label{e:BASE}
		\varrho_n\mb{x}&=\exp\mb{2\pi c\mb{\MNORM{n_M^+}-\MNORM{n_M^-}}}\exp\mb{\iunit 2\pi n^*x}\nonumber\\
		&=\exp\mb{2\pi c\MNORM{n_M^+}}\exp\mb{\iunit 2\pi \mb{n_M^+}^*x_{M^*}^+}\\&\quad\quad\quad\times\exp\mb{-2\pi  c\MNORM{n_M^-}}\exp\mb{\iunit 2\pi \mb{n_M^-}^*x_{M^*}^-}\ft{.}\nonumber
	}
	It is tempting to think of the $\varrho_n$ as basis elements for a tensor product space of a Hardy space on an annulus, with the dual of such a Hardy space. However, we cannot use $\varrho_n$ as such a basis since $n_M^+$ and $n_M^-$ are not independent of each other. Nevertheless, we can decompose $\MC{H}_{A_{M,c}}$ into two generalized Hardy spaces as follows. We define four norms
	\bs{align*}
	{
		\mu_j\mb{f}&\dfn \sup_{y\in A_j}\mb{\int_{\TORUS}\mbm{f\mb{x+\iunit y}}^2\mathrm{d}x}^{\frac 12}\ft{, }f\in L_2\mb{\TORUS}\ft{, }j\in\mbs{1,2,3,4}\ft{, where}\\
		A_1&\dfn\mbs{y\in\VS{\MB{R}}\SR {y_{M^*}^-\in\VS{\mb{-c,c}}\ft{, }y_{M^*}^+\in\VS{\mb{c,\infty}}}}\ft{,}\\
		A_2&\dfn\mbs{y\in\VS{\MB{R}}\SR {y_{M^*}^-\in\VS{\mb{-c,c}}\ft{, }y_{M^*}^+\in\VS{\mb{-\infty,-c}}}}\ft{,}\\
		A_3&\dfn\mbs{y\in\VS{\MB{R}}\SR {y_{M^*}^-\in\VS{\mb{-c,c}}\ft{, }y_{M^*}^+\in{\mb{c,\infty}\times \mb{-\infty,-c}}}}\ft{,}\\
		A_4&\dfn\mbs{y\in\VS{\MB{R}}\SR {y_{M^*}^-\in\VS{\mb{-c,c}}\ft{, }y_{M^*}^+\in{\mb{-\infty,-c}\times \mb{c,\infty}}}}\ft{.}
	}
	Of course, for all $f\in L_2\mb{\TORUS}$ the norms $\mu_j\mb{f}$ are not finite but they are so at least for some Fourier polynomials. The spaces $H_j$, $j\in\mbs{1,2,3,4}$ are the completion \wrt the norms $\mu_j$ above. E.g. using $\mu_1$, it holds for all $f\in H_1$
	\bs{align*}
	{
	\mu_1\mb{f}^2=&\sup_{y\in A_1}\mb{\int_{\TORUS}\mbm{f\mb{x+\iunit y}}^2\mathrm{d}x}=\sup_{y\in A_1}\sum_{n\in\VS{\MB{Z}}}\exp\mb{-4\pi n^*y}\mbm{\varphi_n^*f}^2\\
		=&\sup_{y\in A_1}\sum_{n\in\VS{\MB{Z}}}\exp\mb{-4\pi \mb{n^-_M}^*y_{M^*}^--4\pi \mb{n^+_M}^*y_{M^*}^+}\mbm{\varphi_n^*f}^2\\
		=&\sup_{y^+_{M^*}\in \VS{\mb{c,\infty}}}\sum_{n\in\VS{\MB{Z}}}\exp\mb{4\pi c \mnp[{}]{n^-_M}-4\pi \mb{n^+_M}^*y_{M^*}^+}\mbm{\varphi_n^*f}^2\\
		=&\sum_{\substack{n\in\VS{\MB{Z}}\\n_M^+\in\VS{\mbro{0,\infty}}}}\exp\mb{4\pi c \mnp[{}]{n^-_M}-4\pi c \mnp[{}]{n^+_M}}\mbm{\varphi_n^*f}^2=	\sum_{\substack{n\in\VS{\MB{Z}}\\n_M^+\in\VS{\mbro{0,\infty}}}}\mbm{\varphi_n^*A_{M,c}f}^2\ft{.}	
		}
		Similar calculations for the other three norms show then that the spaces $H_j$, $j\in\mbs{1,2,3,4}$ disjointly partition the space $\MC{H}_{A_{M,c}}$ \wrt the dual coordinate up to $n=0$. Since $E_M^+$ is a one dimensional subspace of $\MB{R}^2$, always two of the spaces contain only the constant functions (note that $n_M^+=0$ implies $n=0$), say, $H_3$ and $H_4$. Then all vectors in the spaces $H_1$ and $H_2$ are holomorphic functions on $\TORUS+\iunit A_1$ and on $\TORUS+\iunit A_2$, respectively.
}

	\section{The Koopman operator is nuclear}
\label{s:TRACEKOOP}
 We set for each $r>0$
 \[\PERT[\MC{T}]\dfn\mbs{T\colon\TORUS\rightarrow\TORUS\SR T\ft{ extends holomorphically and boundedly on }\STRIP{r}}\ft{.}\]
 For every $T\in\PERT[\MC{T}]$ the Koopman operator
 \[\MC{K}_T\colon L_2\mb{\TORUS}\rightarrow L_2\mb{\TORUS}\colon f\mapsto f\circ T\]
 is well-defined by differentiability of $T$. It is well-known that the operator $\MC{K}_T$ acting on $L_2\mb{\TORUS}$ is not compact. We say that two maps $f$, $g\in\PERT[\MC{T}]$ are $C^1$-close if the distance
 \[d\mb{f,g}\dfn\sup_{z\in\STRIP{r}}\mnp[{}]{f(z)-g(z)}+\sup_{z\in\STRIP{r}}{\mnp[{}]{{\DIFF[z]f-\DIFF[z]g}}}\]
 is small. In this section we revisit the proof of Faure and Roy \cite{faure2006ruelle}. They showed that $\MC{K}_T$, acting on the Hilbert space $\MC{H}_{A_{M,c}}$, (see \de{HS}), is nuclear of order $0$ if $T$ is sufficiently $C^1$-close to a hyperbolic matrix $M\in \SL[\MB{Z}]{2}$ for some $c>0$.\newline
 We recall that a linear operator $\MC{L}\colon \MC{H}\rightarrow\MC{H}$ on a Hilbert space $\MC{H}$ with norm $\mnp[\MC{H}]{\cdot}$ is called nuclear of order $0$ if it can be written as a sum $\MC{L}= \sum_{n\in\MB{N}}d_n\psi_{1,n} \psi_{2,n}^*$
 with $\inf\mbs{p>0\SR\sum_{n\in\MB{N}}\mbm{d_n}^p<\infty}=0$ and $\psi_{1,n},\psi_{2,n} \in\MC{H}$, $\mnp[\MC{H}]{\psi_{1,n}},\mnp[\MC{H}]{\psi_{2,n}}\le 1$, $d_n\in\MB{C}$, $n\in\MB{N}$ \cite[II, \S 1, n$^\circ$1, p.4]{Grothen_1966}. In particular, such an operator is trace class, hence bounded and admits a trace $\tr\MC{L}\dfn\sum_{n\in\MB{N}}e_n^*\MC{L}e_n$, invariant for any choice of orthonormal basis $e_n$, $n\in\MB{N}$ of $\MC{H}$. Moreover, one can show that $\tr\MC{L}$ equals the sum, including multiplicity (dimension of corresponding generalized eigenspace), over the spectrum $\ft{sp}\mb{\MC{L}}$ of $\MC{L}$. The Fredholm determinant, defined for small enough $z\in\MB{C}$ by
 \bs{align}
 {
 	\label{e:FD}
 	\det\mb{1-z\MC{L}}\dfn \exp\mb{-\sum_{n=1}^\infty\frac {z^n}n\tr{\MC{L}^n}}\ft{,}
 }
 extends to an entire function in $z$, having zeroes at $z=\lambda^{-1}$, $\lambda\in\ft{sp}\mb{\MC{L}}\setminus\mbs{0}$ of same order as the multiplicity of $\lambda$.
%

\bo{Nuclearity of $\MC{K}_T$}{theorem}
{
	\label{t:TRACECLASS}
	Let \ftBranch{$M\in\SL[\MB{Z}]{2}$}{$d\in\MB{N}$ and $M\in\MAT[\MB{Z}]{d}$} be hyperbolic and let $r>0$. Then there exist constants $\delta_M>0$ and $c_1>0$ such that for each $T\in\PERT[T]$ with $d(T,M)\le\delta_M$ the map
	\[\MC{K}_{T}\colon\MC{H}_{A_{M,c_1}}\rightarrow \MC{H}_{A_{M,c_1}}\colon f\mapsto f\circ T\]
	defines a nuclear operator of order $0$. In particular, there exists $c_2>0$ depending only on $c_1$, $M$, and $\mnp[{}]{\cdot}$ so that for each $n_1$, $n_2\in\VS{\MB{Z}}$
	\[\mbm{\mba{\varrho_{n_1},\MC{K}_{T}\varrho_{n_2}}_{\MC{H}_{A_{M,c_1}}}}\le\exp\mb{-2\pi c_2\mb{\mnp[{}]{n_1}+\mnp[{}]{n_2}}}\ft{.}\]
	Moreover, $c_1<r$, and the constants $c_1$, $c_2$, and $\delta_M$ can be chosen arbitrarily small.
}

For every $n_1$, $n_2\in\VS{\MB{Z}}$, we set 
\bs{equation}
{
	\label{e:OSC}
	\FUNC{n_1}{n_2}\dfn \mba{\varphi_{n_1},\MC{K}_T\varphi_{n_2}}_{L_2\mb{\TORUS}}\ft{.}
}
Estimating this "oscillatory integral" is central for \tm{TRACECLASS}. In the case $T=M$, we have simply
\bs{align}
{
	\label{e:SIMPLY}\FUNC[M]{n_1}{n_2}=\bs{cases}{1&\ft{if }M^*n_2=n_1\\0&\ft{if }M^*n_2\ne n_1}\ft{.}
}
The strategy of the proof is as follows. We get an upper bound for $\mbm{\FUNC{n_1}{n_2}}$ in \lem{ESTIMATE}, taking advantage of the holomorphicity of $T$. 
In \lem{BNDMATEL}, we compare the contribution of $n_1$ and $n_2$ in the expanding and contracting directions, using here essentially the hyperbolicity of $M$. Combining both results, we obtain a weaker bound on $\mbm{\FUNC{n_1}{n_2}}$ in \pp{INTBOUND}, which finally allows for the proof of \tm{TRACECLASS}.\\
For every $n\in\VS{\MB{Z}}$ and $y\in\VS{\MB{R}}$ any solution $x\in\TORUS$ so that
\bs{equation}
{
	\label{e:MW}
	\exp\mb{-2\pi\mb{n^*\DIFF{T}y}}=\int_{\TORUS}\exp\mb{-2\pi\mb{n^*\DIFF[z]{T}y}}\mathrm{d}z
}
is denoted by $x_n\mb{y}$. Since the integrand is continuous in $y$ such a solution exists by the Mean Value Theorem.

\bo{Upper bound on $\mbm{\FUNC{n_1}{n_2}}$ (I)}{lemma}
{
	\label{l:ESTIMATE}
	Let \ftBranch{}{$d\in\MB{N}$, }$r>0$. Then, there exists $C\ge 0$ so that for each $T\in\PERT[T]$ and $n_1$, $n_2\in\VS{\MB{Z}}$ and $y\in\VS{\mb{-r,r}}$, recalling (\ref{e:OSC}), we have
	\[\mbm{\FUNC{n_1}{n_2}}\le \exp\mb{2\pi\mb{-n_2^*\DIFF[x_{n_2}(y)]{T}y+n_1^*y+C{d(T,0)}\mnp[{}]{y}^3\mnp[{}]{n_2}}}\ftn{.}\]
}
\bs{proof}
{
	By definition \[\FUNC{n_1}{n_2}=\mba{\varphi_{n_1},\MC{K}_T\varphi_{n_2}}_{L_2\mb{\TORUS}}=\int_{\TORUS}\exp\mb{\ftn{i}2\pi\mb{n_2^*T(x)-n_1^*x}}\mathrm{d}x\ft{.}\]
	Since $T\in\PERT[T]$, the $\VS{\MB{Z}}$-invariance of the integrand follows. By holomorphicity of $T$ on $\STRIP{r}$, we can change the path of integration to $x\mapsto x+\ftn{i}y$ for every $y\in\VS{\mb{-r,r}}$. Therefore for any $y\in\VS{\mb{-r,r}}$ \[\mbm{\FUNC{n_1}{n_2}}\le \int_{\TORUS}\exp\mb{2\pi\mb{n_1^*y-\Im\mb{n_2^*T(x+\iunit y)}}}\mathrm{d}x\ftn{,}\]
	where $\Im$ is the imaginary part.
	We expand $T$ (or rather its lift to $\VS{\MB{R}}$) at $x\in\TORUS$ in a Taylor series to the second order. This yields
	\[T\mb{x+\iunit y}=T(x)+\iunit \DIFF{T}y+P(x+\iunit y)+R_2\mb{x+\iunit y}\ft{.}\] 
	Here, $P(x+\iunit y)$ is the second order term of the expansion which is $\VS{\MB{R}}$-valued, and $R_2$ is the remainder of the series expansion. We find therefore
	\[\Im T(x+\ftn{i}y)=\DIFF{T}y+\Im R_2 \mb{x+\iunit y}\ft{.}\]
	 Since $T$ is holomorphic we find a constant $C>0$ independent of $T$ such that  \[\mbm{n_2^*R_2\mb{x+\iunit y}}\le C{d(T,0)}\mnp[{}]{n_2}{\mnp[{}]{y}^3}\ft{.}\] We are left with the evaluation of
	\[\int_{\TORUS}\exp\mb{-2\pi\mb{n^*\DIFF[z]{T}y}}\mathrm{d}z\ft{.}\]
	Using (\ref{e:MW}) yields the result.
}

The following abbreviation is used in the remaining section. We set for each $y\in\VS{\MB{R}}$
\bs{equation}
{
	\label{e:CHOICE}
	\mbm{y}_{M}\dfn \ANORM{y_M^+}- \ANORM{y_M^-}\ft{.}
}
\bo{Directional inequality}{lemma}
{
	\label{l:BNDMATEL}
		Let \ftBranch{$M\in\SL[\MB{Z}]{2}$}{$d\in\MB{N}$ and $M\in\MAT[\MB{Z}]{d}$} be hyperbolic. Let
		$\epsilon>0$ and $\kappa\ge 0$ and let $R\colon\VS{\MB{R}}\rightarrow \MB{R}_{\ge 0}$ be a map such that for all $z\in\VS{\MB{R}}$ with $\mnp[{}]{z}<\epsilon$ it holds
		 \[R(z)\le \kappa\mnp[{}]{z}\ft{.}\] 
	Then there exists $c_M>0$ such that for all $\kappa<c_M$ there exist $c_1>c_2>0$ such that for all $n_1$, $n_2\in\VS{\MB{Z}}$ there exists $y_{n_1,n_2}\in\VS{\MB{R}}$ independent of $R$ with $\mnp[{}]{y_{n_1,n_2}}<\epsilon$ such that it holds
	\[-c_1\mb{\mbm{n_1}_M-\mbm{n_2}_M}-\mb{n_2^*M-n_1^*}y_{n_1,n_2}+\mnp[{}]{n_2}R\mb{y_{n_1,n_2}}\le -c_2\mb{\mnp[{}]{n_1}+\mnp[{}]{n_2}}\ft{.}\]
}
\bs{proof}
{
	We assume $0<c_2\le c_1$. For $n_{1}=n_{2}=0$ there is nothing to prove. For every $\mb{y_1,y_2}\in\VS{\MB{R}}$ we set $\mnp{\mb{y_1,y_2}}\dfn\sqrt{y_1^2+y_2^2}$. We let $0<\tilde{c}_1\le 1\le\tilde{c}_2$ such that 
	\bs{align}
	{
		\label{e:EQUIVNORM}
		\tilde{c}_2^{-1}\mb{\MNORM{y^+_M}+\MNORM{y^-_M}}\le\MNORM{y}\le\tilde{c}_1^{-1}\mnp[2]{y}\ft{,}\quad\ft{for all $y\in\VS{\MB{R}}$. }
	}Whenever $n_{2}\ne 0$ we find a linear map $M_a$ such that $M_{a}n_2=M^*n_2-n_1$ and whenever $n_{1}\ne 0$ we find a linear map $M_b$ such that $M_{b}n_1=M^*n_2-n_1$.	For now we let $\tilde{\kappa}>0$ be a variable which will be fixed later on, independently of $n_1$ and $n_2$. We consider the following four cases
	\bd{multicols}{2}
	{
		\bo{label=(\alph*)}{enumerate}
		{
			\item\label{CASEA} $\mnp[{}]{n_2}>0$ and $\mnp{n_2}\ge \mnp{n_1}$
			\bo{label=(\roman*)}{enumerate}
			{\item\label{CASEA:I} $\mnp[]{M_an_2}\ge\tilde{\kappa}\mnp[]{n_2}$\ft{,}\item\label{CASEA:II} $\mnp[]{M_an_2}<\tilde{\kappa}\mnp[]{n_2}$\ft{,}}
			\item\label{CASEB} $\mnp[{}]{n_1}>0$ and $\mnp{n_1}\ge \mnp{n_2}$\bo{label=(\roman*)}{enumerate}{\item\label{CASEB:I} $\mnp[]{M_bn_1}\ge\tilde{\kappa}\mnp[]{n_1}\ft{,}$\item\label{CASEB:II} $\mnp[]{M_bn_1}<\tilde{\kappa}\mnp[]{n_1}$.}
		}
	}
	We assume Case \ref{CASEA}\ref{CASEA:I}. For every $\delta> 0$ we let
	\[y=\delta M_a^{*}\frac{n_2}{\mnp[{}]{n_2}}\ft{.}\]
	
	It follows, using (\ref{e:EQUIVNORM}), that
	\bs{align}
	{
		\label{e:CAI1}
		-\mb{n_2^*M-n_1^*}y
		=-n_2^*M_a^*y
		\le-\tilde{c}_1^2\mnp[{}]{M_a^*{n_2}}\mnp[{}]{y}\ft{.}
	}
	We recall $\mbm{\cdot}_M$ from (\ref{e:CHOICE}). Using that $c_1+c_2>0$ and that \ref{CASEA} holds, we estimate
	\bs{align*}
	{
		-c_1\mb{\mbm{n_1}_M-\mbm{n_2}_M}&\le c_1\mb{\mnp[{}]{n_1}+\mnp[{}]{n_2}}\\
		&=-c_2\mb{\mnp[{}]{n_1}+\mnp[{}]{n_2}}+\mb{c_1+c_2}\mb{\mnp[{}]{n_1}+\mnp[{}]{n_2}}\\
		&\le -c_2\mb{\mnp[{}]{n_1}+\mnp[{}]{n_2}}+2\mb{c_1+c_2}\mnp[{}]{n_2}\ft{.}
	}
	Using \ref{CASEA}\ref{CASEA:I} and the assumed bound on $R$ for $\mnp[{}]{y}<\epsilon$, we have
	\bs{align}
	{
		\label{e:CAI2}
		-c_1\mb{\mbm{n_1}_M-\mbm{n_2}_M}&-\mnp[{}]{n_2}\mb{\tilde{c}_1^2\mnp{M_a^*\frac{n_2}{\mnp[{}]{n_2}}}\mnp[{}]{y}-R\mb{y}}\le\nonumber\\& -c_2\mb{\mnp[{}]{n_1}+\mnp[{}]{n_2}}+\mb{2(c_1+c_2)-\mb{\tilde{c}_1^2\tilde{\kappa}-\kappa}\mnp[{}]{y}}\mnp[{}]{n_2}\ft{.}
	}
	
	We put $c_M\dfn\tilde{c}_1^2\tilde{\kappa}$. Any value $\mnp[{}]{y}\in\mb{0,\epsilon}$ can be attained by controlling $\delta$. Assuming that $c_M>\kappa$, it follows from (\ref{e:CAI1}) and (\ref{e:CAI2}) that
	\bs{align}
	{
		\label{e:C1C2BOUND}
		0<c_1+c_2< \frac{c_M-\kappa}{2}\epsilon\ft{.}
	}
	The reasoning in Case \ref{CASEB}\ref{CASEA:I} is completely analogous and yields the same bounds on $c_1+c_2$.\\
	In Case \ref{CASEA}\ref{CASEA:II} and \ref{CASEB}\ref{CASEA:II}, we take $y=0$, where $R(0)=0$ by assumption on $R$. We assume now Case \ref{CASEA}\ref{CASEA:II}. We find, using (\ref{e:EQUIVNORM}),
    \bs{align}
    {
    	\label{e:NE}
	    \mnp[{}]{\mb{M_a^*n_2}^+_M}+\mnp[{}]{\mb{M_a^*n_2}^-_M}&\le \tilde{c}_2\mnp[{}]{M_a^*n_2}\nonumber\\
	    &<\tilde{c}_2\tilde{\kappa}\mnp[{}]{n_2}\le\tilde{c}_2\tilde{\kappa}\mb{\mnp[{}]{n^+_{2,M}}+\mnp[{}]{n^-_{2,M}}}\ft{.}
	}
    We have
	\[
		\ANORM{\mb{M_a^*{n_2}}^+_M}=\ANORM{M^*n_{2,M}^+-n_{1,M}^+}\quad\ft{and}\quad\ANORM{\mb{M_a^*{n_2}}^-_M}=\ANORM{M^*n_{2,M}^--n_{1,M}^-}\ft{.}
	\]
		
	Recalling (\ref{e:CE}), this allows the estimate
	\bs{align*}
	{
		\ANORM{\mb{M_a^*{n_2}}^+_M}+\ANORM{\mb{M_a^*n_2}^-_M}&\ge \ANORM{M^*n_{2,M}^+}-\ANORM{n_{1,M}^+} -\ANORM{M^*n_{2,M}^-}+ \ANORM{n_{1,M}^-}\\
		&\ge \lambda_M \ANORM{n_{2,M}^+}-\lambda_M^{-1} \ANORM{n_{2,M}^-} - \ANORM{n_{1,M}^+} + \ANORM{n_{1,M}^-}\ft{.}
	}
	Together with (\ref{e:NE}) we find therefore
	\[-\mbm{n_1}_M=-\ANORM{n_{1,M}^+} + \ANORM{n_{1,M}^-}< -\mb{\lambda_M-\tilde{\kappa} \tilde{c}_2}\ANORM{n_{2,M}^+}+\mb{\lambda_M^{-1}+\tilde{\kappa} \tilde{c}_2}\ANORM{n_{2,M}^-}\ft{.}\]
	We set
	\[\kappa_+\dfn\lambda_M-\tilde{\kappa} \tilde{c}_2-1\quad\ft{and}\quad \kappa_-\dfn 1-\lambda_M^{-1}-\tilde{\kappa} \tilde{c}_2\ft{.}\]
	We finally estimate
	\bs{align*}
	{
		-c_1\mb{\mbm{n_1}_M-\mbm{n_2}_M}< -c_1\kappa_+\ANORM{n_{2,M}^+}-c_1\kappa_-\ANORM{n_{2,M}^-}\ft{.}
	}
	Note that we have $\kappa_+>\kappa_-$ because $\lambda_M>1$. Assuming that $c_1\kappa_-\ge 2c_2$, we find
	\bs{align*}
	{
		-c_1\kappa_-\ANORM{n_{2,M}^+}-c_1\kappa_-\ANORM{n_{2,M}^-}< -c_1\kappa_-\mnp[{}]{n_2}\le -2c_2\mnp[{}]{n_2}\le -c_2\mb{\mnp[{}]{n_1}+\mnp[{}]{n_2}}\ft{.}
	}
	In Case \ref{CASEB}\ref{CASEA:II} we consider the bounds
	\bs{align*}
	{
		\mbm{n_2}_M+\lambda_M\ANORM{n_{1,M}^-}-\lambda_M^{-1}\ANORM{n_{1,M}^+}&\le \mnp[{}]{\mb{\mb{M^*}^{-1}M_b^*n_1}_M^+}+\mnp[{}]{\mb{\mb{M^*}^{-1}M_b^*n_1}_M^-}\\
		&\le\tilde{c}_2\mnp[{}]{\mb{M^*}^{-1}M_b^*n_1}<\tilde{\kappa} \tilde{c}_2\mnp[{}]{\mb{M^*}^{-1}}\mnp[{}]{n_1}\ft{.}
	}
	Therefore $\kappa_-$ is replaced by $1-\lambda_M^{-1}-\mnp[{}]{\mb{M^*}^{-1}}\tilde{\kappa} \tilde{c}_2$ which we require to be positive. Since $\mnp[{}]{\mb{M^*}^{-1}}>1$, this yields the stronger conditions
	\bs{align}
	{
		\label{e:C2IN}
		0<\tilde{\kappa}<\frac{1-\lambda_M^{-1}}{\mnp[{}]{\mb{M^*}^{-1}}\tilde{c}_2}\quad\ft{and}\quad c_2\le\frac{1-\lambda_M^{-1}-\mnp[{}]{\mb{M^*}^{-1}}\tilde{\kappa} \tilde{c}_2}2c_1\ft{.}
	}
	Any such choice for $\tilde{\kappa}$ is independent of $n_1$ and $n_2$ and fixes $c_M$. Using (\ref{e:C2IN}) for an upper bound on $c_2$, we find
	\[c_1+c_2\le \frac{3-\lambda_M^{-1}}{2}c_1\ft{.}\]
	Using (\ref{e:C1C2BOUND}), we require the stronger condition
	\[0<c_1< \frac{c_M-\kappa}{{3-\lambda_M^{-1}}}\epsilon\ft{.}\]
	Therefore the choices of $c_1$ and $c_2$ are valid if $\kappa<c_M$. They depend only on $\epsilon$, $M$ and $\mnp[{}]{\cdot}$ and not on $n_1$ or $n_2$.
}

\bo{Upper bound on $\mbm{\FUNC{n_1}{n_2}}$ (II)}{proposition}
{
	\label{p:INTBOUND}
	Let \ftBranch{$M\in\SL[\MB{Z}]{2}$}{$d\in\MB{N}$ and $M\in\MAT[\MB{Z}]{d}$} be hyperbolic and let $r>0$. Then there exist constants $0<\delta_M$ and $0<c_2<c_1$ such that for each $n_1$, $n_2\in\VS{\MB{Z}}$ and each $T\in\PERT[T]$ with $d(T,M)\le\delta_M$ it holds that
	\bs{align*}
	{
		\exp\mb{-2\pi c_1\mb{\mbm{n_1}_M-\mbm{n_2}_M}}\mbm{I_{n_1,n_2}\mb{T}}\le\exp\mb{-2\pi c_2\mb{\mnp[{}]{n_1}+\mnp[{}]{n_2}}}\ft{.}
	}
}
\bs{proof}
{
	By \lem{ESTIMATE} there is a constant $C>0$ independent of $T$ such that for each $y\in \VS{\mb{-r,r}}$ and $n_1$, $n_2\in\VS{\MB{Z}}$ it holds that
	\bs{align}
	{
		\label{e:INEQ}
		\mbm{I_{n_1,n_2}\mb{T}}\le\exp\mb{2\pi\mb{-n_2^*\DIFF[x_{n_2}(y)]{{T}}y+n_1^*y+Cd(T,0)\mnp[{}]{y}^3\mnp[{}]{n_2}}}\ft{\ft{.}}
	}
	We rewrite
	\bs{align*}
	{
		n_2^*\DIFF[x_{n_2}\mb{y}]{T}y&=n_2^*My+n_2^*\DIFF[x_{n_2}(y)]\mb{T-M} y\ft{,}
	}
	and set
	\[R(y)\dfn \bs{cases}{\frac{n_2^*}{\mnp[{}]{n_2}}\DIFF[x_{n_2}(y)]\mb{M-T} y+C{d(T,0)}\mnp[{}]{y}^3 &\ft{if }n_2\ne 0\\ 0&\ft{if }n_2=0}\ft{.}\]
	Let $\delta_M> 0$ and assume that $d(T,M)\le \delta_M$. We choose $0<\epsilon\le r$ sufficiently small such that for all $y\in\VS{\MB{R}}$ with $\mnp[{}]{y}<\epsilon$ there is $\kappa>0$ such that
	\[\mbm{R(y)}\le\kappa\delta_M\ft{.}\]
	Since $d(T,0)\le d(T,M)+d(M,0)\le \delta_M+d(M,0)$ this choice of $\epsilon$ is independent of $T$. \lem{BNDMATEL} applied to $M$ and $\mbm{R}$ gives $c_1$, $c_2$ and $y_{n_1,n_2}\in\VS{\MB{R}}$ for which the \rhs of (\ref{e:INEQ}) fulfills the desired inequality. 
}

\bo{Proof of \tm{TRACECLASS}}{proof}
{
	\pp{INTBOUND} yields $0<\delta_M$ and $0<c_2<c_1$ such that if $d(T,M)\le\delta_M$ it holds
	\bs{equation}
	{
		\label{e:BOUND}
		C_{n_1}C_{n_2}^{-1}\mbm{I_{n_1,n_2}\mb{T}}\le \exp\mb{-2\pi c_2\mb{\mnp[{}]{n_1}+\mnp[{}]{n_2}}}\ft{,}
	}
	where \[C_n\dfn\exp\mb{-2\pi c_1\mb{\MNORM{n_M^+}-\MNORM{n_M^-}}}\ft{,}\quad n\in\VS{\MB{Z}}\ft{.}\]
	
	We put $c\dfn c_1$ and $M$ in Definitions \ref{d:IM} and \ref{d:HS}, giving a linear map $A_{M,c_1}$ and a Hilbert space $\MC{H}_{A_{M,c_1}}$.
 Recalling (\ref{e:FB}), and assuming that $\MC{K}_{T}\colon \MC{H}_{A_{M,c_1}}\rightarrow \MC{H}_{A_{M,c_1}}$ is well-defined, we have
	\bs{align}{\label{e:THM}\mbm{\mba{\varrho_{n_1},\MC{K}_{T}\varrho_{n_2}}_{\MC{H}_{A_{M,c_1}}}}&=\mbm{\mba{\varphi_{n_1},A_{M,c_1}\MC{K}_{T}A_{M,c_1}^{-1}\varphi_{n_2}}_{L_2\mb{\TORUS}}}\nonumber\\
		&=C_{n_1}C_{n_2}^{-1}\mbm{I_{n_1,n_2}\mb{T}}\ft{.}}
	Using (\ref{e:BOUND}) to estimate the \rhs, the bound in \tm{TRACECLASS} follows. We next obtain well-definedness and nuclearity of order $0$ of $\MC{K}_{T}$. Let $f\in\MC{H}_{A_{M,c_1}}$ and put $g\dfn A_{M,c_1}f$. We have then
	\bs{align*}
	{
		\MC{K}_{T}f\in\MC{H}_{A_{M,c_1}}&\Leftrightarrow A_{M,c_1}\MC{K}_{T}f\in L_2\mb{\TORUS}\Leftrightarrow\sum_{n\in\VS{\MB{Z}}}\mbm{\varphi_{n}^*A_{M,c_1}\MC{K}_{T}f}^2<\infty\\
		&\Leftrightarrow\sum_{n_1\in\VS{\MB{Z}}}\mbm{\sum_{n_2\in\VS{\MB{Z}}}\varphi_{n_1}^*A_{M,c_1}\MC{K}_{T}A_{M,c_1}^{-1}\varphi_{n_2}\varphi_{n_2}^*g}^2<\infty\\&\Leftrightarrow \sum_{n_1\in\VS{\MB{Z}}}\mbm{\sum_{n_2\in\VS{\MB{Z}}}C_{n_1}C_{n_2}^{-1}I_{n_1,n_2}\mb{T}\varphi_{n_2}^*g}^2<\infty\ft{.}
	}
	Using (\ref{e:BOUND}) and the Cauchy-Schwartz inequality, it follows that
	\bs{align*}
	{
		\sum_{n_1\in\VS{\MB{Z}}}\mbm{\sum_{n_2\in\VS{\MB{Z}}}C_{n_1}C_{n_2}^{-1}I_{n_1,n_2}\mb{T}\varphi_{n_2}^*g}^2\le \mb{\sum_{n\in\VS{\MB{Z}}}e^{-4\pi c_2{\mnp[{}]{n}}}}^2\mnp[L_2\mb{\TORUS}]{g}^2<\infty\ft{.}
	}
	This gives the well-definedness  of $\MC{K}_{T}$. Now, using the Cauchy-Schwartz inequality, we have
	\[\mbm{\mba{\varrho_n,\MC{K}_Tf}_{\MC{H}_{A_{M,c_1}}}}^2\le\sum_{m\in\VS{\MB{Z}}}\mbm{\mba{\varrho_n,\MC{K}_T\varrho_m}_{\MC{H}_{A_{M,c_1}}}}^2\mnp[\MC{H}_{A_{M,c_1}}]{f}^2\ft{.}\]
	Using (\ref{e:THM}) and (\ref{e:BOUND}) to bound $\mbm{\mba{\varrho_n,\MC{K}_T\varrho_m}_{\MC{H}_{A_{M,c_1}}}}$, we find a constant $C>0$ such that
	\[\mbm{\mba{C\exp\mb{2\pi c_2\mnp[{}]{n}}\varrho_n,\MC{K}_Tf}_{\MC{H}_{A_{M,c_1}}}}\le \mnp[\MC{H}_{A_{M,c_1}}]{f}\ft{.}\] 
	This allows the representation of $\MC{K}_T$ as 
	\[\MC{K}_Tf=\sum_{n\in\VS{\MB{Z}}}C^{-1}\exp\mb{-2\pi c_2\mnp[{}]{n}}\mba{C\exp\mb{2\pi c_2\mnp[{}]{n}}\varrho_n,\MC{K}_Tf}_{\MC{H}_{A_{M,c_1}}}\varrho_n\ft{,}\] 
	from which nuclearity of order $0$ follows. Finally, a brief inspection of the proofs for \lem{BNDMATEL} and \pp{INTBOUND} gives the statement about the constants.
}
	\section{Non-trivial resonances for the Koopman operator}
\label{s:SPECRUELLE}
Given any hyperbolic matrix $M\in\SL[\MB{Z}]{2}$, we find by \tm{TRACECLASS} constants $0<\delta_M$ and $c>0$ such that for each map $T\in\PERT[\MC{T}]$, satisfying $d(T,M)\le\delta_M$, the operator $\MC{K}_T$ acting on the Hilbert space $\MC{H}_{A_{M,c}}$ is nuclear of order $0$. Therefore it has a well-defined trace
\bs{equation}
{
	\label{e:TRACE}
	\tr\MC{K}_{T}\dfn\sum_{n\in\VS{\MB{Z}}}\mba{\varrho_n,\MC{K}_{T}\varrho_n}_{\MC{H}_{A_{M,c}}}\ft{.}
}
The map $T$ is an Anosov diffeomorphism (for all small enough $\delta_M$), by structural stability \cite[Theorem 9.5.8]{hasselblatt_first_2003}. Then the map $T$ has the same number $N_M=\mbm{\det\mb{\id-M}}$ of fixed points as the matrix $M$. We recall a well-known result \cite[Proposition 9]{faure2006ruelle}.

\bo{Trace formula for $\MC{K}_T$}{lemma}
{
	\label{l:TRACEFORMULA}
	Let \ftBranch{$M\in\SL[\MB{Z}]{2}$}{$d\in\MB{N}$ and $M\in\MAT[\MB{Z}]{d}$} be hyperbolic and let $r>0$. Then there exist constants $\delta_M>0$ and $c>0$ such that for each $T\in\PERT[T]$ with $d(T,M)\le\delta_M$, letting $\MC{K}_T$ act on $\MC{H}_{A_{M,c}}$, it holds
	\[\tr\MC{K}_{T}=\sum_{T\mb{x}=x}\mbm{\det\mb{\id-\DIFF{T}}}^{-1}\ft{.}\]
}
For the convenience of the reader, we give a proof:\bs{proof}
{
	 Using \tm{TRACECLASS} gives constants $c>0$ and $\delta_M>0$ and well-definedness of $\MC{K}_{T}$. For small enough $\delta_M>0$, by structural stability and \lem[II]{FP}, the map $\id-T$ can be partitioned into $N_M$ surjective submaps. In particular, there are diffeomorphisms $y_j\colon D_j\rightarrow \TORUS$, $D_j\subseteq \TORUS$, $1\le j\le N_M$ such that $\id-T=\bigsqcup_{j=1}^{N_M}y_j$. Then, using (\ref{e:FB}), we have for each $n\in\VS{\MB{Z}}$
	\bs{align*}
	{
		\mba{\varrho_n,\MC{K}_{T}\varrho_n}_{\MC{H}_{A_{M,c}}}&=
		\mba{\varphi_n,A_{M,c}\MC{K}_{T}A_{M,c}^{-1}\varphi_n}_{L_2\mb{\TORUS}}=\int_{{\TORUS}}\exp\mb{\iunit 2\pi n^*\mb{T-\id}\mb{x}}\mathrm{d}x\\
		&=\sum_{j=1}^{N_M}\int_{y_j^{-1}\mb{\TORUS}}\exp\mb{\iunit 2\pi n^*y_j(x)}\mathrm{d}x\\
		&=\sum_{j=1}^{N_M}\int_{\TORUS}\frac{\exp\mb{\iunit 2\pi n^*z}}{\mbm{\det\mb{\id-\DIFF[y_j^{-1}\mb{z}]{T}}}}\mathrm{d}z\ft{.}
	}	
	For $N\in\MB{N}$ and $z\in\TORUS$ the following sum
	\[D_N\mb{z}\dfn\sum\limits_{\substack{n\in\VS{\MB{Z}}\\\mnp[{}]{z}\le N}}\exp\mb{\iunit 2\pi n^*z}\]
	is the \ftBranch{$2$}{$d$}-dimensional analogue of the Dirichlet kernel \cite[p.13]{katznelson_introduction_2004}. Together with (\ref{e:TRACE}), this yields immediately 
	\[\tr \MC{K}_{T}=\lim_{N\to\infty}\sum\limits_{\substack{n\in\VS{\MB{Z}}\\\mnp[{}]{n}\le N}}\mba{\varrho_n,\MC{K}_{T}\varrho_n}_{\MC{H}_{A_{M,c}}}=\sum_{T\mb{x}=x}\mbm{\det\mb{\id-\DIFF[x]{T}}}^{-1}\ft{.}\]
}

Using \lem{TRACEFORMULA}, and the definitions $(\ref{e:DYNDET})$ and (\ref{e:FD}) for the dynamical determinant and Fredholm determinant, respectively, we see directly that 
\bs{align}
{
	\label{eq:EQUALITY}
	\det\mb{1-z\MC{K}_T}=d_T\mb{z}\ft{.}
}
The Ruelle resonances correspond to the zeroes of the Fredholm determinant, hence to the inverses of the non-zero eigenvalues of $\MC{K}_T$.

\bs{remark}
{
	\label{r:SP}
	In view of \eq{EQUALITY} and the relation of the Ruelle resonances of $T$ to the eigenvalues of $\MC{K}_T$, one may ask how the spectrum of $\MC{K}_T$ would be affected if we let $\MC{K}_T$ act on a different Banach space. The following relates a part of the eigenvalues of two linear operators sharing a common dense subspace and is due to a proof of Baladi and Tsujii \cite[Appendix A]{baladi_dynamical_2008}. Consider two separable Banach spaces $\mb{\MC{B}_1,\mnp[1]{\cdot}}$ and $\mb{\MC{B}_2,\mnp{\cdot}}$. This induces two other Banach spaces
	\[\mb{\MC{B}_1+\MC{B}_2,\mnp[+]{\cdot}}\ft{ and }\mb{\MC{B}_1\cap\MC{B}_2,\mnp[\cap]{\cdot}}\ft{, where}\]
	\bs{align*}
	{
		\mnp[+]{f}&\dfn\inf\mbs{\mnp[1]{f_1}+\mnp[2]{f_2}\SR f_1\in\MC{B}_1\ft{, }f_2\in\MC{B}_2\ft{, }f=f_1+f_2}\ft{ and}\\
		\mnp[\cap]{f}&\dfn\max\mbs{\mnp[1]{f},\mnp[2]{f}}\ft{.}
	}
	Suppose that $\MC{B}_\cap$ is dense in $\MC{B}_1$ and $\MC{B}_2$. Let $\MC{L}\colon \MC{B}_+\rightarrow \MC{B}_+$ be a linear map which preserves the spaces $\MC{B}_\cap$, $\MC{B}_1$ and $\MC{B}_2$ and is a bounded linear operator on the restrictions $\MC{L}_{|\MC{B}_1}$ and $\MC{L}_{|\MC{B}_2}$. Then the part of the eigenvalues of $\MC{L}_{|\MC{B}_1}$ and $\MC{L}_{|\MC{B}_2}$ coincide which lies outside the closed disc with radius larger to both essential spectral radii. Moreover, the corresponding generalized eigenspaces of $\MC{L}_{|\MC{B}_1}$ and $\MC{L}_{|\MC{B}_2}$ coincide and are contained in $\MC{B}_\cap$.\\
	For the applications that we have in mind, the map $\MC{L}$ is just the Koopman or transfer operator, defined on $\MC{B}_1$ and $\MC{B}_2$, respectively, extended to the space $\MC{B}_+$.
}

The spectrum $\ft{sp}\mb{\MC{K}_T}$ of $\MC{K}_T$ on $\MC{H}_{A_{M,c}}$ is invariant under complex conjugation since $T$ is real. The constant functions on $\TORUS$\ftBranch{}{, $d\in\MB{N}$} are all fixed by $\MC{K}_T$. Therefore we have $1\in \ft{sp}\mb{\MC{K}_T}$. If we take $T=M^k$, $k\in\MB{N}$ in \lem{TRACEFORMULA}, it follows that $\tr \MC{K}_T=1$. Hence, the dynamical determinant is just $d_T\mb{z}=1-z$, also noted in \cite[p.3]{Rugh_1996}. We find immediately that $1$ is the only Ruelle resonance. We show now that this finding is non-generic in the following sense. The rest of this section is devoted to an idea of Naud \cite{naud_2015}. We put for every $r>0$
\bs{align}
{\label{e:BS}\PERT\dfn\mbs{T\in\PERT[\MC{T}]\SR \ft{The lift of }T\ft{ to }\VS{\MB{R}}\ft{ is }\VS{\MB{Z}}\ft{-periodic}}\ft{.}}
Endowed with the uniform norm this is a Banach space.
\bo{Non-trivial Ruelle resonances (I)}{theorem}
{
	\label{t:NTRS}
	\ftBranch{Let $M\in\SL[\MB{Z}]{2}$}{Let $d\in\MB{N}$ and $M\in\MAT[\MB{Z}]{d}$}
	 be hyperbolic. For each $r>0$ there exists an open and dense set $\MC{G}\subseteq\PERT$ such that
	 the linear functional
	 \[ B_M\colon\PERT\rightarrow\MB{R}\colon \psi\mapsto{N_M}^{-1}\sum_{Mx=x}\tr\mb{\mb{\id-M}^{-1}\DIFF[x]\psi}\]
	 never vanishes on $\MC{G}$. For all $\psi\in\MC{G}$ there exists $\epsilon_0>0$ such that for all $\epsilon<\epsilon_0$
	\[\tr\MC{K}_{M+\epsilon \psi}=1+\epsilon B_M\mb{\psi}+O\mb{\epsilon^2}\ft{.}\]
	In particular,  	\[\spec\mb{\MC{K}_{M+\epsilon\psi}}\setminus\mbs{0,1}\ne \emptyset\ft{.}\]
}

\bo{Real analyticity of fixed points}{lemma}
{
	\label{l:CONT}
	\ftBranch{Let $M\in\SL[\MB{Z}]{2}$}{Let $d\in\MB{N}$ and $M\in\MAT[\MB{Z}]{d}$} be hyperbolic and $r>0$. Then for all $\psi\in\PERT$ the fixed points of the map
	\[M+\delta\psi\]
	are real analytic functions of $\delta$ where $\delta$ lies in a real neighborhood of $0$.
}
\bs{proof}
{
	We set for $\delta\in\MB{R}$
	\[F\mb{\delta,x}\dfn Mx+\delta\psi(x)-x\ft{.}\]
	 We fix a point $y_j\dfn\mb{0,x_j}$ where $x_j$, $1\le j\le N_M$ is a fixed point of $M$.
	 By construction, the map $F$ has a holomorphic extension to $\MB{C}\times\STRIP{r}$. Since $M$ is hyperbolic, we have $\det\DIFF[x_j]\mb{F\mb{\delta,\cdot}}\ne 0$ for small $\delta$. We apply the Holomorphic Implicit Function Theorem \cite[Theorem 1.4.11]{krantz_function_2001} on $F$ with $F\mb{y_j}=0$. This yields a holomorphic function $x_j\mb{\delta}$ with $x_j(0)=x_j$. It is obviously real analytic for $\delta\in\MB{R}$ in a neighborhood of $0$ and $x_j(0)=x_j$.
}

\bo{Proof of \tm{NTRS}}{proof}
{
	Let $\delta\in\MB{R}$ and $\psi\in\PERT$ and set $\tilde{M}\dfn M+\delta\psi$. We choose $\delta$ small in \lem{CONT} which gives for each fixed point $x$ of $M$ a real analytic function $\tilde{x}$ with $\tilde{x}\mb{0}=x$. Using a Taylor expansion on $\tilde{x}$ at $0$, we have
	\[\tilde{x}\mb{\delta}=x+O(\delta)\ft{.}\]
	Using real analyticity of the derivative $\DIFF{\psi}$, we have
	\[{\DIFF[x]\psi-\DIFF[\tilde{x}\mb{\delta}]\psi}= O(\delta)\ft{.}\]
	We write now for each fixed point $x$ of $M$
	\bs{align*}
	{
		\mbm{\det\mb{\id-\DIFF[\tilde{x}\mb{\delta}]{\tilde{M}}}}&=\mbm{\det\mb{\id-M-\delta\DIFF[x]{\psi}+\delta\mb{\DIFF[x]{\psi}-\DIFF[\tilde{x}\mb{\delta}]{\psi}}}}\\
		&=N_M\mbm{\det\mb{\id-\mb{\id-M}^{-1}\mb{\delta\DIFF[x]{\psi}+\mb{\delta\DIFF[x]{\psi}-\delta\DIFF[\tilde{x}\mb{\delta}]{\psi}}}}}\\
		&=N_M\mbm{\det\mb{\id-\delta\mb{\id-M}^{-1}\DIFF[x]{\psi}+O\mb{\delta^2}}}\\
		&=N_M\mb{1-\delta\tr\mb{\mb{\id-M}^{-1}\DIFF[x]{\psi}}+O\mb{\delta^2}}\ft{.}
	}
	We have by \lem{TRACEFORMULA} for $\delta$ small enough
	\[\tr \MC{K}_{\tilde{M}}=1+\frac{\delta}{N_M}\sum_{Mx=x}\tr\mb{\mb{\id-M}^{-1}\DIFF[x]{\psi}}+O\mb{\delta^2}\ft{.}\]
	Now we set
	\[B_M\colon\PERT\rightarrow\MB{R}\colon \psi\mapsto{N_M}^{-1}\sum_{Mx=x}\tr\mb{\mb{\id-M}^{-1}\DIFF[x]\psi}\ft{.}\]
	We next check that this is a non-trivial linear functional. Note that formally $\psi\mb{\id-M}=2$. However, no non-zero linear map is in the space of additive perturbations $\MC{B}_r$. We denote by $v_j$, \ftBranch{$j\in\mbs{1,2}$}{$1\le j\le d$} the $j$-th column of the matrix $\mb{\mb{\id-M}^*}^{-1}$ and we fix now $j$\ftBranch{}{ such that $v^*_j\ne 0$}. Let \ftBranch{$\psi_0\colon \TORUS[1]+\iunit\mb{-r,r}\rightarrow\MB{C}$ be holomorphic and bounded}{$\psi\in \MC{B}_{1,r}$}. For every $\mb{x_1,x_2}\eqqcolon x\in\TORUS$ we put
	\bs{equation*}
	{
		\psi\mb{x}\dfn\psi_0\mb{x_j}v_j\ft{.}
	}
	By construction, we have $\psi\in\PERT$ and we evaluate 
	\bs{equation*}
	{
	 	B_M\mb{\psi}=\frac{v_j^*v_j}{N_M}\sum_{Mx=x}\psi_0^{(1)}\mb{x_j}\ft{.}
	}
	 The \rhs is a finite sum and by taking for $\psi_0$ a suitable Fourier polynomial (e.g. a shifted sinus with sufficiently high frequency), we can establish $B_M\mb{\psi}\ne 0$.
	 We set $\MC{G}\dfn B_M^{-1}\mb{\MB{R}\setminus\mbs{0}}$. By continuity of $B_M$, the set $\MC{G}$ is open and dense in $\PERT$. 
}
	\section{Non-trivial resonances for the transfer operator}
\label{s:SPECTRANSFER}
 As before, we consider maps $T\in\PERT[T]$, $r>0$ which are sufficiently $C^1$-close to a hyperbolic linear map $M\in\ftBranch{\SL{2}}{\MAT{d}\ft{, }d\in\MB{N}}$.
 We turn to the adjoint of $\MC{K}_{T}$, acting on the dual space $\MC{H}_{A_{M,c}}^*$, which we denote by $\MC{L}_{T}$.

\bo{Transfer operator}{lemma}
{
	\label{l:TO}
	Let \ftBranch{$M\in\SL[\MB{Z}]{2}$}{$d\in\MB{N}$ and $M\in\MAT[\MB{Z}]{d}$} be hyperbolic and let $r>0$. Then there exist constants $0<\delta_M$ and $c>0$ such that for each $T\in\PERT[T]$ with $d(T,M)\le\delta_M$ the map
	\[\MC{L}_{T}\colon\MC{H}_{A_{M,c}}^*\rightarrow\MC{H}_{A_{M,c}}^*\colon f\mapsto{
		\ftBranch
		{{\frac{f}{\mbm{\det\DIFF[{}]{T}}}\circ T^{-1}}}
		{\sum_{j=1}^{\mbm{\det{M}}}{\frac{f^*}{\mbm{\det\DIFF[{}]{T}}}\circ T_j^{-1}}}}\]
	defines a nuclear operator of order $0$, conjugate to $\MC{K}_{T}$. In particular, \[\ft{sp}\mb{\MC{L}_{T}}=\ft{sp}\mb{\MC{K}_{T}}\ft{.}\]
}
\bs{proof}
{
	By \tm{TRACECLASS} there is $0<\delta_M$, $c>0$ and $\MC{H}_{A_{M,c}}$ such that $\MC{K}_{T}$ acting on $\MC{H}_{A_{M,c}}$ is nuclear of order $0$ if $d(T,M)\le\delta_M$. The same can be said about its adjoint, acting on $\MC{H}_{A_{M,c}}^*$ (e.g. see \cite[p. 77]{ryan_introduction_2002}). The trace of $\MC{K}_T$ and $\MC{L}_T$ coincide, so does their Fredholm determinant, and hence their resonances. By definition of the adjoint, $\forall f^*\in \MC{H}_{A_{M,c}}^*$, $\forall g\in\MC{H}_{A_{M,c}}\colon \mb{\MC{L}_{T}f}^*\mb{g}=f^*\mb{\MC{K}_{T}g}$. Using \lem{DS}, it holds
\ftBranch
{\bs{align*}
{
			f^*\mb{\MC{K}_Tg}&=\mba{A_{M,c}^{-2}f,\MC{K}_{T}g}_{\MC{H}_{A_{M,c}}}=\int_{\TORUS}\mb{A_{M,c}^{-1}\bar{f}}\mb{x}\mb{A_{M,c}\MC{K}_{T}g}\mb{x}\mathrm{d}x\\
			&=\int_{\TORUS}\bar{f}\mb{x}\mb{\MC{K}_{T}g}\mb{x}\mathrm{d}x=\int_{\TORUS}\mb{\frac{\bar{f}}{\mbm{\det\DIFF[{}]T}}\circ T^{-1}}\mb{x}g\mb{x}\mathrm{d}x\\
			&=\mba{A_{M,c}^{-2}\mb{\frac{f}{\mbm{\det\DIFF[{}]T}}\circ T^{-1}},g}_{\MC{H}_{A_{M,c}}}=\mb{\frac{f}{\mbm{\det\DIFF[{}]T}}\circ T^{-1}}^*\mb{g}\ft{.}}}
{\bs{align*}
	{
		\mba{f,\MC{K}_{T}g}_{\MC{H}_{A_{M,c}}}&=\mba{A_{M,c}f,A_{M,c}\MC{K}_{T}g}_{L_2\mb{\TORUS}}=\int_{\TORUS}\mb{A_{M,c}\bar{f}}\mb{x}\mb{A_{M,c}\MC{K}_{T}}g\mb{x}\mathrm{d}x\\
		&=\int_{\TORUS}\mb{A_{M,c}^2\bar{f}}\mb{x}\mb{\MC{K}_{T}g}\mb{x}\mathrm{d}x\\
		&=\sum_{j=1}^{N_{\id+T}}\int_{T_j^{-1}\mb{\TORUS}}\mb{A_{M,c}^2\bar{f}}\mb{x}\mb{\MC{K}_{T_j}g}\mb{x}\mathrm{d}x\\
		&=\sum_{j=1}^{N_{\id+T}}\int_{\TORUS}\frac{A_{M,c}^2\bar{f}}{\mbm{\det\DIFF[{}]T_j}}\circ T_j^{-1}\mb{x}g\mb{x}\mathrm{d}x\\
		&=\mba{\sum_{j=1}^{N_{\id+T}}A_{M,c}^{-2}\mb{\frac{A_{M,c}^2f}{\mbm{\det\DIFF[{}]T_j}}\circ T_j^{-1}},g}_{\MC{H}_{A_{M,c}}}\ft{.}
}}\ftBranch{}{The claim follows $\mbm{\det\DIFF[{}]{T_j}}\circ T_j^{-1} = \mbm{\det\DIFF[{}]{T}}\circ T_j^{-1}$ and $N_{\id+M}=\mbm{\det M}$.} 
}

By \lem{TO}, recalling (\ref{e:FB}), and \lem{TRACEFORMULA} it holds
\[\tr \MC{L}_T=\sum_{n\in\VS{\MB{Z}}}\MC{L}_T\varrho_n^*\mb{\varrho_n}=\sum_{T(x)=x}\mbm{\det\mb{\id-\DIFF{T}}}^{-1}\ft{.}\]
We have the equality
\[d_T\mb{z}=\det\mb{1-z\MC{K}_T}=\det\mb{1-z\MC{L}_T}\ft{.}\]
We give now analogously to \tm{NTRS} a spectral result for the transfer operator (recall $\PERT$ from (\ref{e:BS})).

\bo{Non-trivial Ruelle resonances (II)}{lemma}
{
	\label{l:NTRSTO}
	\ftBranch{Let $M\in\SL[\MB{Z}]{2}$}{Let $d\in\MB{N}$ and $M\in\MAT[\MB{Z}]{d}$} be hyperbolic. For each $r>0$ there exists an open and dense set $\MC{G}\subseteq\PERT$ such that for all $\psi\in\MC{G}$ there exists $\epsilon_0>0$ such that for all $\epsilon\le\epsilon_0$
	\[\spec\mb{\MC{L}_{M+\epsilon\psi}}\setminus\mbs{0,1}\ne\emptyset\ft{.}\]
}
\bs{proof}
{
	By \tm{NTRS} we know that under every perturbation $\psi\in\MC{G}$ there is $\epsilon_0>0$ such that we find for all $\epsilon\le\epsilon_0$ non-trivial Ruelle resonances. Using \lem{TO} for well-definedness of $\MC{L}_{M+\epsilon\psi}$ and for the relation $\ft{sp}\mb{\MC{L}_{T}}=\ft{sp}\mb{\MC{K}_{T}}$, the result follows.
}
Clearly, the Lebesgue measure (by \rem{DS}, the constant density $1$) is fixed by $\MC{L}_{M}$. This does not persist under a generic perturbation of $M$. However, the spectral relation in \lem{TO} implies that $\MC{L}_T$ fixes some functionals in $\MC{H}_{A_{M,c}}^*$. In particular, using \rem{SP}, we can apply \cite[Theorem 3]{blank_ruelleperronfrobenius_2002} to our transfer operators $\MC{L}_M$ and $\MC{L}_T$. Hence, the eigenvalue $1$ of $\MC{L}_T$ is simple and the projector $\Pi_1^*$ onto the corresponding eigenspace of $\MC{L}_T$ gives us the SRB measure
\[\mu_{\ft{SRB}}\dfn {{\Pi_1^*}1^*}\ft{,}\]
in the usual sense. (It is absolutely continuous \wrt Lebesgue measure in the unstable direction.)\\
We finish this section by showing the existence of non-zero perturbations $\psi\in\PERT$ which allow the determinant $\det\mb{M+\epsilon\DIFF[x]{\psi}}$ to remain constant or to vary along $x\in\TORUS$.
%
\bo{Volume under perturbations}{lemma}
{
	\label{l:PERTENDO}
	Let \ftBranch{}{$d\in\MB{N}$ and} $r>0$ and let \ftBranch{$M\in\SL[\MB{Z}]{2}$}{$M\in\MAT[\MB{Z}]{d}$} be hyperbolic. Then there exist non-zero maps $\psi\in\PERT$ in each following case:
	\bo{label=(\roman*)}{enumerate}
	{
		\item\label{l:PERTENDO:I} For all $\epsilon>0$ and all $x\in\TORUS$ it holds $\det\mb{M+\epsilon \DIFF\psi}=1$. 
		\item\label{l:PERTENDO:II} For all $\epsilon>0$ and Lebesgue almost all $x\in\TORUS$ it holds $\mbm{\det\mb{M+\epsilon \DIFF\psi}}\ne 1$. 
	}
	In particular, the map $\psi$ can be chosen such that for all small $\epsilon>0$ the corresponding transfer operator
	\[\MC{L}_{M+\epsilon\psi}\]
	admits non-trivial Ruelle resonances.
}
\bs{proof}
{
	We prove first Claim \ref{l:PERTENDO:I}, including the statement about the non-trivial Ruelle resonances. We will apply \lem[I]{DETPRES}. We choose $\ftBranch{j\in\mbs{1,2}}{1\le j\le d}$, $r>0$ and let $\phi\colon \TORUS[1]+\iunit\mb{-r,r}\rightarrow\MB{C}$ be a holomorphic and bounded map.
	For $\alpha\in\VS{\MB{R}}$ we set for every $\mb{x_1,x_2}\eqqcolon x\in\TORUS$ \[\psi_{\phi,\alpha}\mb{x}\dfn\mb{\alpha_1\phi\mb{x_j}, \alpha_2\phi\mb{x_j}}\ft{.}\]
	We put \ftBranch{$d\dfn 2$}{$d$}, $j$, $T\dfn M$, $\phi$ and $T_\phi\dfn \psi_{\phi,\alpha}$ (e.g. as lift to $\VS{\MB{R}}$) in \lem{DETPRES}.
	Since $M$ is a constant matrix, say, $M=\bs{pmatrix}{a&b\\c&d}$ for suitable $a,b,c,d\in\MB{Z}$, we can write Condition \ref{l:DETPRES} \ref{l:DETPRES:I} as
	\ftBranch
	{
		\bs{equation}{\label{e:ALPHA}\alpha_1 d=\alpha_2 b\quad\ft{if }j=1\quad\ft{or}\quad\alpha_1 c=\alpha_2 a\quad\ft{if }j=2\ft{.}}
	}
	{
		\[\sum_{i=1}^d\mb{-1}^{i}\alpha_i\det\square_{i,j}\mb{M}=0\ft{.}\]
	}Hence,\ftBranch{}{ for $d>1$,} we have non-zero solutions in $\alpha$ independent of $x$. We choose such a solution $\alpha$ and take $\psi=\psi_{\phi,\alpha}$. Then $\psi\in\PERT$ which yields $\det\mb{M+\epsilon \DIFF\psi}=1$ for every $\epsilon>0$.
	We are free to choose any suitable $\phi$. In particular, \tm{NTRS} yields a linear functional $B_M$
	and a dense subset $\MC{G}\subseteq\PERT$ on which $B_M$ is non-zero. We have to make sure that $\psi\in\MC{G}$. Then for $\epsilon$ small $\MC{L}_{M+\epsilon\psi}$ admits non-trivial Ruelle resonances by \lem{NTRSTO}. To this end, we evaluate $B_M$ at $\psi$ which yields 
		 \[ B_M\mb{\psi}=B_M\mb{\psi_{\phi,\alpha}}={N_M}^{-1}\sum_{Mx=x}\tr\mb{\mb{\id-M}^{-1}\DIFF[x]\psi_{\phi,\alpha}}	 	=\frac{v_j^*\alpha}{N_M}\sum_{Mx=x}\phi^{(1)}\mb{x_j}\ft{,}
	\]
	where $v_j^*$ is the $j$-th row of $\mb{\id-M}^{-1}$. The sum over the fixed points of $M$ can be made non-zero by a suitable Fourier polynomial. Now we have
	\[
	v_1^*\alpha = \frac{{\mb{1-d}\alpha_1+c\alpha_2}}{\det\mb{\id-M}}\quad\ft{or}\quad v_2^*\alpha = \frac{{b\alpha_1+\mb{1-a}\alpha_2}}{\det\mb{\id-M}}\ft{.}
	\]
	Using (\ref{e:ALPHA}), we find
	\[
	v_1^*\alpha = \frac{{\mb{c-b+\frac bd}\alpha_2}}{\det\mb{\id-M}}\quad\ft{or}\quad v_2^*\alpha = \frac{{\mb{b-c+\frac ca}\alpha_1}}{\det\mb{\id-M}}\ft{.}
	\]
	Both equations can never be zero since $M$ is not diagonal. We prove now Claim \ref{l:PERTENDO:II} by modifying the map $\psi$. For $\delta\in\MB{R}\setminus\mbs{0}$ we set $\tilde{\alpha}\dfn\alpha+\delta w_j$, where $w_j$ is the $j$-th column of $M$ and put $\widetilde{\psi}\dfn \psi_{\phi,\tilde{\alpha}}$. We have
	\[\det\mb{M+\epsilon\DIFF[x]\widetilde{\psi}}=\det\mb{M+\epsilon \DIFF[x]\psi+\epsilon\DIFF[x]\mb{\widetilde{\psi}-\psi}}=1+\delta\epsilon\phi^{(1)}\mb{x_j}\ft{.}\]
	Since $\phi$ is not constant, the \rhs differs from $1$ (and $-1$) for Lebesgue almost all $x$. Since $v_j^*\tilde{\alpha}=v_j^*\alpha+\delta v_j^*w_j\ne 0$ for the right choice of the sign of $\delta$, we have $B_M\mb{\widetilde{\psi}}\ne 0$.
	
}
\appendix
\section*{Appendix}
\setcounter{section}{1}
\setcounter{theorem}{0}
For the readers convenience we give a proof of a well-known result:
\bo{Fixed points}{lemma}
{
	\label{l:FP}
	Let \ftBranch{$M$ be $2\times 2$ integer matrix}{$d\in\MB{N}$ and $M\in\MAT[\MB{Z}]{d}$} acting on $\TORUS$. Assume that $\det\mb{\id-M}\ne 0$. Then the following holds:
	\bo{label=(\roman*)}{enumerate}
	{
		\item\label{l:FP:I} The number $N_M$ of fixed points of $M$ is given by $N_M=\mbm{\det\mb{\id-M}}$.
		\item\label{l:FP:II} There exists a disjoint partition $D_j\subseteq\TORUS$, $1\le j\le N_M$ of $\TORUS$ such that the maps $y_j\colon D_j\rightarrow \TORUS\colon x\mapsto \mb{\id-M}x$ are bijections.
	}	
}
\bs{proof}
{
	We let $\id-M$ act on the cover $\VS{\MB{R}}$. The linear map $\id-M$ sends a fundamental region of $\TORUS$, e.g. $\VS{\mbro{0,1}}$, to a convex polytope having a non-zero volume given by $\mbm{\det\mb{\id-M}}$. Each fixed point of $M$ on $\TORUS$ is mapped by $\id-M$ to an element of $\VS{\MB{Z}}$, and the number of integer points contained in the polytope is just given by its volume. Claim \ref{l:FP:I} follows.\\
	Let $v_1$, $v_2\in\VS{\MB{Z}}$ be two different such integer points in the polytope. Now assume that there are $f_1$, $f_2\in\VS{\mbro{0,1}}$ such that
	\[\mb{\id-M}^{-1}\mb{f_1-f_2}\equiv\mb{\id-M}^{-1}\mb{v_1-v_2}\quad\mb{\ft{mod } \VS{\mbro{0,1}}}\ft{.}\]
	The \rhs is mapped to a fixed point of $M$ on $\TORUS$, implying that $f_1-f_2$ is an integer point, which is only possible if $f_1=f_2$. Therefore, $v_1=v_2$, which contradicts the assumption, and Claim \ref{l:FP:II} follows.
}

For $d\in\MB{N}$ and every real $d\times d$ matrix $M$ we denote by $\square_{i,j}\mb{M}$, $1\le i$, $j\le d$ the submatrix arising by removing the $i$-th row and $j$-th column from $M$. 

\bo{Determinant preserving transformation}{lemma}
{
	\label{l:DETPRES}
	Let $d\in\MB{N}$, and let $T\colon\MB{R}^d\rightarrow\MB{R}^d$ and $\phi\colon \MB{R}\rightarrow \MB{R}$ be differentiable maps. Fix $1\le j\le d$ and $\alpha\in\MB{R}^d$ and set \[T_\phi\colon\MB{R}^d\rightarrow\MB{R}^d\colon x\mapsto\mb{\alpha_i\phi\mb{x_j}\SR 1\le i\le d }\ft{.}\]
	Then for $x\in\MB{R}^d$ it holds \[\det\DIFF\mb{T+T_\phi}-\det\DIFF\mb{T}=0\]
	if and only if at least one of the conditions holds:
	\bo{label=(\roman*)}{enumerate}
	{
		\item\label{l:DETPRES:I} $\sum_{i=1}^d\mb{-1}^{i}\alpha_i\det\square_{i,j}\mb{\DIFF T}=0$ or
		\item\label{l:DETPRES:II}$\phi^{(1)}\mb{x_j}=0$\ft{.}
	}
}
\bs{proof}
{
	We develop the determinant of $\DIFF\mb{T+T_\phi}$ \wrt the $j$-th column. Since $T_\phi$ depends only on $x_j$ this gives
	\[\det\DIFF\mb{T+T_\phi}=(-1)^j\sum_{i=1}^d (-1)^i\partial_j\mb{T+T_\phi}_i\mb{x}\det\square_{i,j}\mb{\DIFF T}\ft{.}\]
	Hence, it holds
	\bs{align*}
	{
		\det\DIFF\mb{T+T_\phi}-\det\mb{\DIFF T}&=(-1)^j\sum_{i=1}^d (-1)^i\det\square_{i,j}\mb{\DIFF T}\partial_j\mb{T_\phi}_i\mb{x}\\
		&=(-1)^j\phi^{(1)}\mb{x_j}\sum_{i=1}^d\mb{-1}^{i}\alpha_i\det\square_{i,j}\mb{\DIFF T}\ft{.}
	}
	One reads up Claim \ref{l:DETPRES:I} and \ref{l:DETPRES:II} directly from the \rhs.
}

\section*{References}
\printbibliography[heading=none]
}